\def\versiondate{Nov.~8 2000} %

\input math.macros
\input Ref.macros
\input EPSfig.macros

\checkdefinedreferencetrue
\theoremcountingtrue
\sectionnumberstrue
\figuresectionnumberstrue
\forwardreferencetrue
\tocgenerationtrue
\citationgenerationtrue
\hyperstrue
\initialeqmacro


\let\hat\widehat

\def\msnote#1{\ifproofmode {\bf[** #1 **]}\else\relax\fi}
\long\def\comment#1{}
\let\pname=\procname

\def\itemrm#1{\item{{\rm #1}}}
\def\beginitems{\begingroup\parindent=25pt}
\def\enditems{\vskip0pt\endgroup\noindent}

\def\newbottom#1#2#3{{\eightpoint\parindent=0pt\parskip=2pt\footnote{}
 {{\it 2000 Mathematics Subject Classification.}
 #1}\footnote{}{{\it Key words and phrases.} #2}\footnote{}{#3}}}

\def\verts{{\ss V}}
\def\faces{{\ss F}}
\def\edges{{\ss E}}
\def\Dvert{D_\verts}
\def\Dface{D_\faces}
\def\Dedge{D_\edges}
\def\setminu{-}
\def\eps{\epsilon}
\def\Cal#1{{\cal #1}}
\def\le{\leqslant}
\def\leq{\leqslant}
\def\ge{\geqslant}
\def\geq{\geqslant}
\def\H{{\Bbb H}}
\def\hh{{\H^2}}
\def\KK{K}
\def\hbd{\partial\hh}
\def\dist{{\rm dist}}
\def\Aut{{\rm Aut}}
\def\isomh{{\rm Isom(\hh)}}
\def\area{{\rm area}}
\def\d#1{#1^\dagger}    
\def\ev#1{{\Cal #1}}
\def\bd{\partial}
\def\st{:\,}
\def\operatorname#1{\mathop{{\rm #1}}}
\def\cite#1{\ref b.#1/}

\def\BLPSgip{\cite{BLPSgip}}
\def\BLPScrit{\cite{BLPSdeath}}
\def\BLPSusf{\cite{BLPS:usf}}
\def\BLSpert{\cite{BLS:pert}}

\def\assumptions{Let $G$ be a transitive, nonamenable,
planar graph with one end}

\def\assumeBern{Let $T$ be a vertex-transitive
tiling of $\hh$ with finite sided faces, let
$G$ be the graph of $T$, and let $\omega$
be Bernoulli percolation on $G$}

\def\assumeTo{Let $T$ be a vertex-transitive
tiling of $\hh$ with finite sided faces, let
$G$ be the graph of $T$, and let $\omega$
be an invariant percolation on $G$}

\ifproofmode \relax \else\head{} {\versiondate}\fi
\vglue20pt

\title{Percolation in the Hyperbolic Plane}

\author{Itai Benjamini, and Oded Schramm}

\abstract{Following is a study of percolation in the hyperbolic
plane $\hh$ and on regular tilings in the hyperbolic plane.
The processes discussed include Bernoulli site and bond percolation
on planar hyperbolic graphs, invariant dependent percolations on 
such graphs, and Poisson-Voronoi-Bernoulli percolation.
We prove the existence of three distinct nonempty phases for the Bernoulli
processes.
In the first phase, $p\in(0, p_c]$, there are no unbounded clusters, but there is
a unique infinite cluster for the dual process.  In the second phase,
$p\in(p_c,p_u)$, there are infinitely many unbounded clusters for the process
and for the dual process.  In the third phase, $p\in [p_u,1)$,
there is a unique unbounded cluster, and all the clusters of the dual process 
are bounded. 
We also study the dependence
of $p_c$ in the Poisson-Voronoi-Bernoulli percolation
process on the intensity of the underlying Poisson process.}

\newbottom{Primary
82B43
. %
Secondary 
60K35
, %
60D05
.} 
{Poisson process, Voronoi, isoperimetric constant, nonamenable, planarity, Euler formula,
Gauss-Bonnet theorem.}
{Research partially supported by the Sam and Ayala Zacks Professorial Chair (Schramm,
the Weizmann Institute).}

\bsection {Introduction}{s.intro}

The purpose of this paper is to study percolation in the
hyperbolic plane and in transitive planar graphs that are
quasi-isometric to the hyperbolic plane.  

There are several sources where the reader may consult
for background on percolation on $\Z^d$ \cite{Grimmett:percolation}
and $\R^d$ \cite{MR:contperc} and for background on percolation
on more general graphs
\cite{pyond}, \cite{Lyons:percsurvey}, \cite{pyond:recent}.
For this reason, we will be quite brief here.
Background on hyperbolic geometry may be found
in \cite{CFKW} and the references sited there.

\bigskip\noindent{\bf Percolation on planar hyperbolic graphs.}

A graph $G$ is {\bf transitive} if the automorphism
group of $G$ acts transitively on the vertices $\verts(G)$.
An {\bf invariant percolation} on a graph $G$ is
a probability measure on the space of subgraphs of $G$,
which is invariant under the automorphisms of $G$. 
The connected components
of the random subgraph are often called {\bf clusters}.
Of special interest are the {\bf Bernoulli} site and bond
percolation.  The Bernoulli($p$)  bond percolation
is the random subgraph with vertices $\verts(G)$,
and where each edge is in the percolation subgraph
with probability $p$, independently.
In Bernoulli($p$) site percolation, each vertex is
in the percolation subgraph with probability $p$,
independently, and an edge appears in the percolation
subgraph iff its endpoints are present.

Let $G$ be an infinite, connected, planar, transitive graph,
with finite vertex degree.
Each such graph is quasi-isometric with one and only one
of the following spaces: $\Z$, the $3$-regular tree,
the Euclidean plane $\R^2$, and the hyperbolic plane $\hh$
\cite{Babai:growth}.
Each of these classes has its distinct geometry, and random
processes behave similarly on graphs within each class.
On trees, Bernoulli percolation is quite simple, but
there are some interesting results concerning
invariant percolation \cite{Hag:deptree}.
Bernoulli percolation on $\Z^2$ and planar lattices
in $\R^2$ has an extensive theory. (See \cite{Grimmett:percolation}.)
Burton and Keane \cite{Burton-Keane:planar} also obtained
a theory of invariant percolation in $\Z^2$.

A transitive graph $G$ has one end if for every
finite set of vertices $V_0\subset\verts(G)$ there
is precisely one infinite connected component of $G\setminus V_0$.
A transitive planar graph $G$ with one end is quasi-isometric
to $\R^2$ or to $\hh$ (see, e.g., \cite{Babai:growth}
or the proof of \ref p.unim/.(b), below).
Amenability is a simple geometric
property which distinguishes these two possibilities.  $G$
is {\bf amenable} if for every $\eps>0$ there is 
a finite  set of vertices $V_0$ such that $|\partial V_0| <\eps |V_0|$.
A graph quasi-isometric with $\hh$ must be nonamenable.
Grimmett and Newman \cite{GN:treeZ} have shown
that for some parameter values, Bernoulli percolation
on the cartesian
product of $\Z$ with a regular tree of sufficiently high
degree has infinitely many infinite components.  This
cannot happen on amenable graphs, by the
Burton-Keane argument \cite{BK:uni}.

Let $p_u=p_u(G)$ be the infimum 
of the set of $p\in[0,1]$ such that Bernoulli($p$) percolation
on $G$ has a unique infinite cluster a.s.
The critical parameter $p_c=p_c(G)$ is defined
as the infimum of the set of $p\in[0,1]$ such that
Bernoulli($p$) percolation on $G$ has an infinite cluster a.s.
It has been conjectured \cite{pyond} that 
$p_c(G)<p_u(G)$ for
every nonamenable transitive graph $G$.
It has been recently proven by Pak and Smirnova-Nagnibeda 
\cite{PakSN:uniq} that every nonamenable group has some Cayley
graph $G$ satisfying $p_c(G)<p_u(G)$.
Lalley \cite{Lalley:Fuch} proved $p_c(G)<p_u(G)$ for planar Cayley graphs of
Fuchsian groups with sufficiently high genus.
We now show

\procl t.nun
\assumptions.
Then $0<p_c(G) < p_u(G) <1$, for Bernoulli bond or site percolation on $G$.
\endprocl

The hyperbolic plane seems to be a very good testing ground
for conjectures about nonamenable graphs.  The planarity
and hyperbolic geometry help to settle questions that may
be more difficult in general.

\procl t.pu
\assumptions.
Then Bernoulli($p_u$) percolation on $G$ has a unique
infinite cluster a.s.
\endprocl

This contrasts with a result of Schonmann \cite{S:nuni},
which shows that a.s.\ the number of infinite components of Bernoulli($p_u$)
percolation on $T\times \Z$ is either $0$ or $\infty$ when $T$ is a regular tree. 
(If $T$ is a regular tree of sufficiently high degree,
then $p_u(T\times\Z)>p_c(T\times\Z)$, by
\cite{GN:treeZ}, and hence there are infinitely many
infinite components at $p_u$ on $T\times\Z$.)
Peres \cite{P:nonamen} has generalized
this result of \cite{S:nuni} (with a different proof) to nonamenable products.

It is known that on a transitive graph $G$ when $p>p_u(G)$,
Bernoulli($p$) percolation
a.s.\ has a unique infinite cluster.  This has been proven
by H\"aggstr\"om and Peres~\cite{HP:unimon}
in the unimodular setting, and in full generality 
by Schonmann~\cite{Scho:unimon}.  
The proof of \ref t.pu/ will also prove this \lq\lq uniqueness
monotonicity\rq\rq\ in this restricted setting.

The situation at $p_c$ is also known.
The following theorem is from \BLPSgip; see also \BLPScrit.

\procl t.critperc
Let $G$ be a nonamenable graph with a vertex-transitive
unimodular autormorphism group.
Then a.s.\ critical Bernoulli bond or site percolation on $G$ has
no infinite components.
\endprocl
 
Note that  a planar transitive graph has a unimodular automorphism
group (\ref p.unim/).  Hence the above theorem applies to the graphs
under consideration here.

\bigskip\noindent{\bf Percolation in $\hh$.}

Though percolation is usually studied on graphs, and many interesting
phenomena already appear in the graph setup, there are some special
properties that only show in the continuous setting.

We therefore consider the Poisson-Voronoi-Bernoulli percolation model
in the hyperbolic plane $\hh$.  There are two parameters
needed to specify the model, $\lambda>0$ and $p\in[0,1]$.
Given such a pair $(p,\lambda)$, consider a Poisson point
process with intensity $\lambda$ in the hyperbolic
plane.  Such a process gives rise to a Voronoi tessallation
of $\hh$.  Each tile of this tessallation is colored
white [respectively, black] with probability $p$ [respectively, $1-p$], independently.
The union of all the white [resp.\ black] tiles is denoted by
$\hat W$ [resp.\ $\hat B$].
Let $\Cal W$ [resp.\ $\Cal B$] denote the set of pairs
$(p,\lambda)$ such that $\hat W$ [resp.\ $\hat B$] has an 
unbounded component a.s.  For $\lambda>0$,
let 
$$
p_c(\lambda):=\inf \bigl\{p\in[0,1]\st (p,\lambda)\in\Cal W\bigr\}\,.
$$

We prove the following

\procl t.phases
Consider $(p,\lambda)$-Voronoi percolation in $\hh$.
\beginitems
\itemrm{(a)}
If $(p,\lambda)\in\Cal W\cap \Cal B$, then there are a.s.\
infinitely may unbounded components of $\hat W$ and there
are infinitely many unbounded components of $\hat B$.
\itemrm{(b)}
If $(p,\lambda)\in\Cal W\setminus\Cal B$, then a.s.\ there
is a unique unbounded component of $\hat W$, and no unbounded
components of $\hat B$.
\itemrm{(c)}
If $(p,\lambda)\in\Cal B\setminus\Cal W$, then a.s.\ there
is a unique unbounded component of $\hat B$, and no unbounded
components of $\hat W$.
\enditems
\endprocl

Note that, by symmetry,
$$
\Cal B =\bigl\{(p,\lambda)\st (1-p,\lambda)\in\Cal W\bigr\}
\,.
$$

\midinsert
\centerline{\epsfysize=2.4in\epsfbox{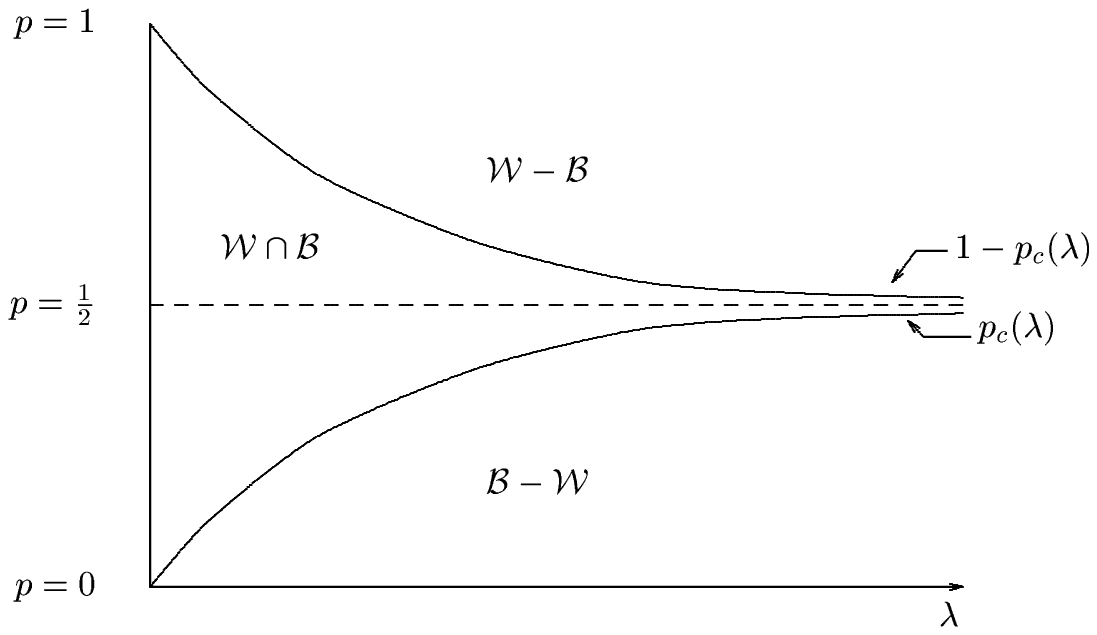}}
\caption{\figlabel{phase}\enspace
The phase diagram of Voronoi percolation.}
\endinsert

\procl t.hber
\beginitems
\itemrm{(a)} 
$ 
0 < p_c(\lambda) \leq
{1\over 2} - {1\over 4\lambda\pi +2}<{1\over 2}
$.
\itemrm{(b)} $\lim_{\lambda\to 0}p_c(\lambda)=0$.
\itemrm{(c)} $p_c(\lambda)$ is continuous.
\itemrm{(d)} $\bigl(p_c(\lambda),\lambda\bigr)\notin\Cal W$ for all $\lambda>0$,
and hence ${\Cal W} =\bigl\{(p,\lambda)\st p>p_c(\lambda)\bigr\}$.
\enditems
\endprocl

We conjecture that $p_c(\lambda)\to 1/2$ as $\lambda\to\infty$.
Note that percolation with parameters $(p,\lambda)$ on $\hh$
is equivalent to percolation with parameters $(p,1)$ on $\hh$ with
the metric rescaled by $\sqrt{\lambda}$.  Hence,
taking $\lambda\to \infty$ amounts to the same as letting the curvature
tend to zero.  This means that Voronoi percolation on $\R^2$
can be seen as a limit of Voronoi percolation on $\hh$.
See \ref q.euc/ for further discussion of this issue.

We generalize the Mass Transport Principle
(\cite{Hag:deptree}, \BLPSgip) to the hyperbolic
plane, and use it as a tool for our investigations.
One consequence of the Mass Transport Principle that we derive
and use is a generalization of Euler's formula
$|V|-|E|+|F|=2$ relating the number of vertices, edges and faces
in a (finite) tiling of the sphere to random (infinite) tilings
of the hyperbolic plane with invariant law.

A collection of open problems is presented at the end of the paper.
Most of these are related to the Voronoi percolation model.
We believe that this process deserves further study.

\bigskip
\noindent {\bf Acknowledgement.}\enspace
We wish to express gratitude to
Harry Kesten, Russ Lyons, Bojan Mohar, Igor Pak, Yuval Peres,
Roberto Schonmann and Bill Thurston
for fruitful conversations and useful advice.

\bsection{Terminology and Preliminaries}{s.prel}

All the graphs that we shall consider in this paper are
locally finite; that is, each vertex has finitely
many incident edges.
The vertices of a graph $G$ will be denoted by $\verts(G)$,
and the edges by $\edges(G)$.

Given a graph $G$, let $\Aut(G)$ denote the group of
automorphisms of $G$.  $G$ is {\bf transitive}
if $\Aut(G)$ acts transitively on the vertices
$\verts(G)$.  $G$ is {\bf quasi-transitive}
if $\verts(G)/\Aut(G)$ is finite; that is,
there are finitely many $\Aut(G)$ orbits
in $\verts(G)$.
A graph $G$ is {\bf unimodular} if $\Aut(G)$ is a unimodular
group (which means that the left-invariant Haar measure
is also right-invariant).  Cayley graphs are unimodular,
and any graph such that $\Aut(G)$ is discrete is unimodular.
See \BLPSgip{} for a further discussion of unimodularity and
its relevance to percolation.

An {\bf invariant percolation} on $G$ is a probability
measure on the space of subgraphs of $G$,
which is $\Aut(G)$-invariant.
A {\bf cluster} is a connected component of the percolation subgraph.

Let $X=\R^2$ or $X=\hh$.
We say that an embedded graph $G\subset X$ in $X$ is {\bf properly} embedded 
if every compact subset of $X$ contains finitely many vertices
of $G$ and intersects finitely many edges.
Suppose that $G$ is an infinite connected graph with
one end, properly embedded in $X$.
Let $\d G$ denote the dual graph of $G$.
We assume that $\d G$ is embedded in $X$ in the standard
way relative to $G$; that is, every vertex $\d v$ of $\d G$
lies in the corresponding face of $G$, and every edge $e\in\edges(G)$
intersects only the dual edge $\d e \in \edges(\d G)$,
and only in one point.
If $\omega$ is a subset of the edges $\edges(G)$,
then $\d \omega$ will denote the set
$$
\d\omega := \bigl\{\d e\st e\notin\omega\bigr\}\,.
$$

Given $p\in[0,1]$ and a graph $G$, we often denote the
percolation graph of Bernoulli($p$) bond percolation on
$G$ by $\omega_p$.

\procl p.unim
\assumptions,
and let $\Gamma$ be the group of automorphism of $G$. 
\beginitems
\itemrm{(a)} $\Gamma$ is discrete (and hence unimodular).
\itemrm{(b)} $G$ can be embedded as a graph $G'$ in the hyperbolic plane
$\hh$ in such a way that the action of $\Gamma$ on $G'$
extends to an isometric action on $\hh$.
Moreover, the embedding can be chosen in such a way
that the edges of $G'$ are hyperbolic line segments.
\enditems
\endprocl

A sketch of the proof of (b) appears in \cite{Babai:growth}.
We include a proof here, for completeness.

\proofof p.unim
By \cite{Mader} or by \cite{Watkins} it follows that 
$G$ is $3$-vertex connected; that is, every finite nonempty set of vertices
$V_0\subset\verts(G)$, $V_0\ne\emptyset$, neighbors with at
least $3$ vertices in $\verts(G)\setminus V_0$.
Therefore, by the extension of Imrich to Whitney's Theorem
\cite{Imrich}, the embedding of $G$ in the plane is unique,
in the sense that in any two embeddings of $G$ in the plane,
the cyclic orientation of the edges going out of the vertices
is either identical for all the vertices, or reversed for all
the vertices.  This implies that an automorphism of $G$
that fixes a vertex and all its neighbors is the identity,
and therefore $\Aut(G)$ is discrete.   For a discrete group,
the counting measure is Haar measure, and is both left- and right-invariant.
Hence $\Aut(G)$ is unimodular.  This proves part (a).

Think of $G$ as embedded in the plane.
Call a component of $S^2-G$ a {\bf face} if its boundary
consists of finitely many edges in $G$.  In each face $f$
put a new vertex $v_f$, and connect it by edges to the vertices
on the boundary of $f$.  If this is done appropriately, then
the resulting graph $\hat G$ is still embedded in the plane.
Note that $\hat G$ together with all its faces forms
a triangulation $T$ of a simply connected domain in $S^2$.
To prove (b) it is enough to produce a triangulation $T'$ of
$\hh$ isomorphic with $T$ such that the elements of $\Aut(T')$ 
extend to isometries of $\hh$ and the edges of $T'$ are hyperbolic
line segments.  There are various ways to do that;
one of them is with circle packing theory.
See, for example, \cite{BeSt}, \cite{HS:hyp}, or
\cite{Babai:growth}.
\Qed

\bsection{The Number of Components}{s.numcomp}

\procl t.nums
\assumptions, and let $\omega$ be an invariant
bond percolation on $G$.
Let $k$ be the number of infinite components of $\omega$,
and $\d k$ be the number of infinite components of $\d \omega$.
Then a.s.
$$
(k,\d k)
\in \Big\{(1,0),(0,1),(1,\infty),(\infty,1),(\infty,\infty)\Big\}
\,.
$$
\endprocl

\procl r.happen
Each of these possibilities can happen.
The case $(k,\d k)=(1,\infty)$ appears when
$\omega$ is the free spanning forest of $G$,
while $(\infty,1)$ is the situation for the
wired spanning forest.
See \BLPSusf.
The other possibilities occur for Bernoulli percolation, as we shall see.
\endprocl

\procl l.nobothfin
\assumptions.  Let $\omega$ be an invariant
percolation on $G$.  If $\omega$ has only finite components
a.s., then $\d\omega$ has infinite components a.s.
\endprocl

The proof will use a result from \BLPSgip, which says
that when the expected degree 
$\E \deg_\omega v$ of a vertex $v$ in an invariant
percolation on a unimodular nonamenable graph $G$ is
sufficiently close to $\deg_G v$, there are infinite
clusters in $\omega$ with positive probability.
(The case of trees was established earlier
in \cite{Hag:deptree}.)

\proof
Suppose that both $\omega$ and $\d\omega$ have
only finite components a.s.
Then a.s.\ given a component $\KK$ of $\omega$, there is a unique
component $\KK'$ of $\d\omega$ that surrounds it.
Similarly, for every component $\KK$ of $\d\omega$,
there is a unique component $\KK'$ of $\omega$ that surrounds it.
Let $\Cal\KK_0$ denote the set of all components of $\omega$.
Inductively, set 
$$
\Cal\KK_{j+1}:= \{ \KK''\st\KK\in\Cal\KK_j\}
\,.
$$
For $\KK\in\Cal\KK_0$ let $r(\KK):=\sup\{j\st\KK\in\Cal\KK_j\}$
be the {\bf rank} of $\KK$, and define $r(v):=r(\KK)$ if
$\KK$ is the component of $v$ in $\omega$.
For each $r$ let $\omega^r$ be the set of edges in $\edges(G)$ incident
with vertices $v\in\verts(G)$ with $r(v)\leq r$.
Then $\omega^r$ is an invariant bond percolation and
$$
\lim_{r\to\infty} \E[\deg_{\omega^r}v] = \deg_G v
\,.
$$
Consequently, by the above result from \BLPSgip, we find that
$\omega^r$ has with positive probability infinite components
for all sufficiently large $r$.
This contradicts the assumption that $\omega$ and
$\d \omega$ have only finite components a.s.
\Qed

The following has been proven in \BLPSgip{}
and \BLPScrit{} in the transitive case.
The extension to the quasi-transitive case
is straightforward.
\msnote{check appropriate ref.  }

\procl t.unnocrit
Let $G$ be a nonamenable quasi-transitive unimodular graph,
and let $\omega$ be an invariant percolation on $G$ which
has a single component a.s.
Then $p_c(\omega)<1$ a.s.
\endprocl

The following has been proven in \BLPSgip.

\procl l.num
Let $G$ be a quasi-transitive nonamenable planar
graph with one end, and let $\omega$ be an invariant
percolation on $G$.  Then a.s.\ the number of infinite
components of $\omega$ is $0,1$ or $\infty$.
\endprocl

For the sake of completeness, we present a (somewhat different) proof here.

\proof
In order to reach a contradiction, assume that with
positive probability $\omega$ has a finite number $k>1$
of infinite components, and condition on that event.
Select at random, uniformly, a pair of distinct infinite components
$\omega_1,\omega_2$ of $\omega$.
Let $\omega_1^c$ be the subgraph of $G$ spanned by the
vertices outside of $\omega_1$, and
let $\tau$ be the set of edges of $G$ that connect vertices
in $\omega_1$ to vertices in the component of $\omega_1^c$
containing $\omega_2$.
Set
$$
\d\tau := \{\d e\st e\in\tau\}
\,.
$$
Then $\d\tau$ is an invariant bond percolation in the dual graph $\d G$.
Using planarity, it is easy to verify that $\d\tau$ is a.s.\
a bi-infinite path.
This contradicts \ref t.unnocrit/, and thereby completes the proof.
\Qed

\procl c.onenou
\assumptions.  Let $\omega$ be an invariant percolation
on $G$.  Suppose that both $\omega$ and $\d\omega$ have
infinite components a.s.
Then a.s.\ at least one among $\omega$ and $\d\omega$
has infinitely many infinite components.
\endprocl

\proof
Draw $G$ and $\d G$ in the plane in such a way that
every edge $e$ intersects $\d e$ in one point, $v_e$
and there are no other intersections of $G$ and $\d G$.
This defines a new graph $\hat G$, whose vertices
are $\verts(G)\cup\verts(\d G)\cup\Big\{v_e\st e\in\edges(G)\Big\}$.
Note that $\hat G$ is quasi-transitive.

Set
$$
\hat\omega :=
\Big\{[v,v_e]\in\edges(\hat G) \st v\in\verts(G),\ e\in\omega\Big\}
\cup
\Big\{[\d v,v_e]\in\edges(\hat G) \st \d v\in\verts(\d G),\ e\notin\omega\Big\}
\,.
$$
Then $\hat\omega$ is an invariant percolation on $\hat G$.
Note that the number of infinite components of $\hat\omega$
is the number of infinite components of $\omega$ plus
the number of infinite components of $\d\omega$.
By \ref l.num/ applied to $\hat\omega$, we find
that $\hat\omega$ has infinitely many infinite components.
\Qed

\proofof t.nums
Each of $k,\d k$  is in $\{0,1,\infty\}$ by \ref l.num/.
The case $(k,\d k)=(0,0)$ is ruled out by \ref l.nobothfin/.
Since every two infinite components of $\omega$ must
be separated by some component of $\d\omega$,
the situation $(k,\d k)=(\infty,0)$ is impossible.
The same reasoning shows that $(k,\d k)=(0,\infty)$
cannot happen.  The case $(k,\d k)=(1,1)$ is ruled out by \ref c.onenou/.
\Qed

\procl t.numsBern
\assumptions, and let $\omega$ be Bernoulli($p$)
bond percolation on $G$.
Let $k$ be the number of infinite components of $\omega$,
and $\d k$ be the number of infinite components of $\d\omega$.
Then a.s.
$$
(k,\d k)
\in \Big\{(1,0),(0,1),(\infty,\infty)\Big\}
\,.
$$
\endprocl

\proof By \ref t.nums/, it is enough to rule out the
cases $(1,\infty)$ and $(\infty,1)$.
Let $K$ be a finite connected subgraph of $G$.
If $K$ intersects two distinct infinite components of
$\omega$, then
$\d\omega\setminu\{\d e\st e\in\edges(K)\}$
has more than one infinite component.
If $k>1$ with positive probability,
then there is some finite subgraph $K$ such that
$K$ intersects two infinite components of $\omega$ with positive
probability.  Therefore, we find that $\d k>1$ 
with positive probability (since the distribution of
$\d\omega\setminu\{\d e\st e\in\edges(K)\}$
is absolutely continuous to the distribution of $\d\omega$).
By ergodicity, this gives $\d k >1$ a.s.
An entirely dual argument shows that $k>1$ a.s.\ when
$\d k>1$ with positive probability.
\Qed

\procl t.BernDual
\assumptions.
Then $p_c(\d G)+p_u(G)=1$ for Bernoulli bond percolation.
\endprocl

\proof
Let $\omega_p$ be Bernoulli($p$) bond percolation on $G$.
Then $\d\omega_p$ is Bernoulli($1-p$) bond percolation on $\d G$.
It follows from \ref t.numsBern/ that the number of infinite
components $\d k$ of $\d \omega$ is
$1$ when $p<p_c(G)$, $\infty$ when $p\in\big(p_c(G),p_u(G)\big)$
and $0$ when $p>p_u(G)$.
\Qed

\proofof t.nun
We start with the proof for bond percolation.
The easy inequality $p_c(G)\geq 1/(d-1)$, where $d$ is the maximal degree of the
vertices in $G$, is well known.

Set $p_c=p_c(G)$.
By \ref t.critperc/, $\omega_{p_c}$ 
has only finite components a.s.
By \ref t.numsBern/, $\d{(\omega_{p_c})}$ has a unique
infinite component a.s.
Consequently, by \ref t.critperc/ again,
$\d{(\omega_{p_c})}$ is supercritical; that is,
$p_c(\d G)< 1-p_c(G)$.  
An appeal to \ref t.BernDual/ now establishes the inequality
$p_c(G)<p_u(G)$.

Since $p_u(G)=1-p_c(\d G)\leq 1 - 1/(\d d-1)$,
where $\d d$ is the maximal degree of the
vertices in $\d G$, we get $p_u(G)<1$, and the proof for bond percolation is complete.

If $\omega$ is site percolation on $G$, let $\omega^b$ be the set of
edges of $G$ with both endpoints in $\omega$.  Then $\omega^b$ is
a bond percolation on $G$.  In this way, results for bond percolation
can be adapted to site percolation.  However, even if $\omega$ is Bernoulli, 
$\omega^b$ is not.  Still, it is easy to check that the above
proof applies also to $\omega^b$.  The details are left to the reader.
\Qed

\proofof t.pu
By \ref p.unim/, $\Aut(G)$ is discrete and unimodular.
By \ref t.BernDual/, $\d{(\omega_{p_u})}$ is critical
Bernoulli bond percolation on $\d G$.
Hence, by \ref t.critperc/, $\d{(\omega_{p_u})}$ has
a.s.\ no infinite components.
Therefore, it follows from \ref t.nums/ that
$\omega_{p_u}$ has a single infinite component.
\Qed

\bsection{Geometric Consequences}{s.s}

We now investigate briefly the geometry of percolation clusters
on vertex-transitive tilings of $\hh$.  Recall that the
{\bf ideal boundary} $\hbd$ of $\hh$ is homeomorphic to the circle $S^1$.
Given a point $o\in\hh$, $\hbd$ can be identified with the space of
infinite geodesic rays starting from $o$.
Let $z_n$ be a sequence in $\hh$.
We say that $z_n$ converges to a point $z$ in $\hbd$,
if the geodesic segments $[o,z_n]$ converge to the ray corresponding to $z$.
(That is, the length of $[o,z_n]$ tends to infinity and the angle at $o$
between the ray from $o$ corresponding to $z$ and the segment $[o,z_n]$
tends to zero.)
One can show that the convergence of $z_n$ does not depend
on the choice of $o$.

\procl t.blim
\assumeBern.
Almost surely, every infinite component of $\omega$ contains
a path that has a unique limit point in the ideal boundary
of $\hh$.
\endprocl

In \BLSpert{} it is shown that for any unimodular transitive graph $G$
and every $p\in[0,1]$, a.s.\ every infinite component of Bernoulli($p$)
percolation on $G$ is transient and simple random walk on it has positive speed.

\proof
Suppose that $\omega$ has infinite components with positive probability.
Let $X(t)$ be a simple random walk on an infinite component of $\omega$.
Then, by the above result of \BLSpert{}, a.s.\ 
there is a $\Lambda>0$ such that
$$
\dist(X(t),X(0))\geq t\Lambda
\,,
\label e.sp
$$
for all sufficiently large $t$,
where $\dist$ is the distance in the hyperbolic metric.
Clearly, we also have,
$$
\dist(X(t),X(s))\leq L|t-s|
\,,
\label e.Lip
$$
for some constant $L$.
These inequalities are enough to conclude that
$X(t)$ tends to a limit in the ideal boundary.
Indeed, fix a polar coordinate system $(r,\theta)$,
where $r=r(p)$ is the hyperbolic distance
$\dist(X(0),p)$ from $p$ to $X(0)$
and $\theta$ is the angle between the segment
$[X(0),p]$ and some fixed ray starting at $X(0)$.
Note that there are constants $c<1<C$
such that for points $p,q\in\hh$
$$
d(\theta(p),\theta(q)) \leq C c^{ r(p) - \dist(p,q)}
\,,
$$
where the distance $d(\theta(p),\theta(q))$ refers to 
the arclength distance on the unit circle.
It therefore follows immediately from \ref e.sp/
and \ref e.Lip/ that
$\{\theta(X(t))\}$ is a Cauchy sequence; that is,
$\lim_t\theta(X(t))$ exists.
This implies that
$X(t)$ tends to a limit in the ideal boundary.
\Qed

There is also a proof that does not use
speed, but instead uses the easier result from
\BLSpert{} that the infinite components are transient.
We now give that proof.

Recall that a metric $d_0$ on
the vertices of a graph $H$ is proper, if every ball of finite
radius contains finitely many vertices.  

\procl l.extl 
Let $G$ be an infinite connected (locally finite) graph.
Then $G$ is transient iff every proper metric on $G$ satisfies 
$$
\sum_{[v,u]\in\edges(G)} d_0(v,u)^2 =\infty
\,.
$$
\endprocl

This easy observation is certainly not new.

\proof
A function $f:\verts(G)\to\R$ is {\bf proper}
if $f^{-1}(K)$ is finite for every compact $K$.
The {\bf Dirichlet energy} of $f$ is
$\Cal D (f):=\sum_{[v,u]\in\edges(G)} \bigl(f(u)-f(v)\bigr)^2$.
It is well known that $G$ is recurrent iff there is a proper function
$f:\verts(G)\to[0,\infty)$ with finite Dirichlet energy.
(See \cite{DoyleSnell} or \cite{Lyons:book}.)

Suppose that $f$ is proper and has finite Dirichlet energy.
Without loss of generality, we may assume that $f(v)\neq f(u)$
if $v\neq u$.  Then define $d_0(v,u):=|f(v)-f(u)|$.
In the other direction, given a proper metric $d_0$,
set $f(u)=d_0(o,u)$, where $o\in\verts(G)$ is arbitrary.
By the triangle inequality it follows that $\Cal D(f)<\infty$.
\Qed

\medbreak\noindent{\it Second proof of \ref t.blim/.\enspace}
We first prove that every transient connected subgraph
$H\subset G$ has a path with a limit in the ideal boundary.
Indeed, consider the Poincar\'e 
(disk) model for the hyperbolic plane $\hh$.
Given two points, $x,y\in\hh$, let $d_E(x,y)$
be the {\it Euclidean\/} distance from $x$ to $y$.
Let $d_E^H$ be the maximal metric on
$\verts(H)$ such that $d_E^H(v,u)\leq d_E(v,u)$ 
whenever $[v,u]\in\edges(H)$ (in other words,
$d_E^H(v,u)=\inf \sum_{j=1}^n d_E(v_{j-1},v_j)$,
where the infimum is taken with respect to all finite
paths $v_0,v_1,\dots,v_n$ in $H$ from $v$ to $u$).

We now show that
$$
\sum_{[v,u]\in\edges(G)} d_E(v,u)^2<\infty
\,.
\label e.admis
$$
Let $o\in\verts(G)$ be a basepoint, and for every
neighbor $u$ of $o$ let $D_u$ be a compact disk
in $\hh$ whose boundary contains the points $o$ and $u$.
Let $P$ be the collection of all disks of the form
$\gamma D_u$ where $u$ neighbors with $o$ and
$\gamma\in\Gamma$ is an automorphism of $G$ acting
on $\hh$ as an isometry.  Since $\Gamma$ acts
discretely on $\hh$,
there is a finite upper bound $M$ for the number of
disks in $P$ that contain any point $z\in \hh$.
Since each disk in $P$ is contained in $\hh$,
it follows that the sum of the Euclidean areas
of the disks in $P$ is finite.  As the square
of the Euclidean diameter of a disk is linearly
related to the Euclidean area, this proves
\ref e.admis/.

Because $H$ is transient, \ref l.extl/ implies that
the metric $d_E^H$ is not proper.  Hence there is an infinite
simple path in $H$ with finite $d_E$-length.
\ref t.blim/ follows, because a.s.\ the infinite components of $\omega$ are transient, by \BLSpert{}.
\Qed

\procl l.dense
\assumeTo.
Let $Z$ be the set of points $z$ in the ideal boundary
$\hbd$ such that there is a path in $\omega$ with limit
$z$.  Then a.s.\  $Z=\emptyset$ or $Z$ is dense in $\hbd$.
\endprocl
\proof
Given a vertex $v\in\verts(G)$, we may consider a polar coordinate
system with $v$ as the origin, and with respect to
this coordinate system $\hbd$ can be thought of as a metric circle of
circumference $2\pi$.  Let $d_v$ denote this metric
of $\hbd$, and let $a(v)$ be the length of the largest component
of $\hbd\setminus Z$, with respect $d_v$.
Note that for vertices $v\in\verts(G)$, the law
of the random variable $a(v)$ does not depend on $v$.
Let $o\in\verts(G)$, let $\epsilon\in(0,1)$,
and let $\delta$ be the probability
that $\eps<a(o)<2\pi-\eps$.

Suppose that $\delta>0$.
Let $R>0$ be very large, and let $x$ be a random-uniform
point on the circle of radius $R$ about $o$ in $\hh$.
Let $v_x$ be the vertex of $G$ closest to $x$.
On the event $\eps<a(o)<2\pi-\eps$,
there is probability greater than $\eps/(4\pi)$
that the geodesic ray from $o$ containing $x$ hits $\hbd$
at a point $x'$ with $d_o(x',Z)\ge \eps/4$; this happens
when $x'$ is within the inner half of the largest
arc of $\hbd\setminus Z$ with respect to $d_o$.
On that event, if $R$ is very large (as a function of
$\eps$), we have $a(v_x)$ as close as we wish to
$2\pi$, and $a(v_x)\ne2\pi$.  However, since
$$
\P\bigl[a(v_x)\in(2\pi-t,2\pi)\bigl]=
\P\bigl[a(o)\in(2\pi-t,2\pi)\bigl]{\to}0\,,\qquad\hbox{as }t\searrow 0\,,
$$
it follows that $\P[a(o)\in (\eps,2\pi-\eps)]=0$.
Consequently, $a(o)\in\{0,2\pi\}$ a.s.

It remains to deal with the case that $Z$ is a single point with
positive probability.  But this is impossible, since the law of
this point would be a finite positive measure on
$\hbd$ which is invariant under  the automorphisms of
$T$, and such a measure does not exist.
(To verify this, observe that such a measure can have no atoms,
since the orbit of any point in $\hbd$ is infinite.
Given any finite positive atomless measure 
on $\hh$ and a small $\eps>0$, there is a bounded
set of points $x\in\hh$ with the property that
for every half plane which contains $x$, the
arc of $\hbd$ associated with it has measure at least $\eps$.
However, the orbit of every point $x\in\hh$ under
the automorphisms of $T$ is infinite.)
\Qed

\procl c.half
\assumeBern.
If $\omega$ has infinite components a.s., then for every
half space $W\subset\hh$, a.s.\ $\omega\cap W$ has infinite
components.
\endprocl

The significance of percolation in hyperbolic half-spaces
was noted by \cite{Lalley:Fuch}.

\proof
This follows immediately from
\ref t.blim/ and \ref l.dense/. \Qed


\procl c.connectivity \procname{Connectivity decay}
Let $T$ be a vertex-transitive
tiling of $\hh$ with finite sided faces, let
$G$ be the graph of $T$,
let $p<p_u(G)$, and let $\omega=\omega_p$
be Bernoulli($p$) percolation on $G$.
Let $\tau(v,u)$ denote the probability
that $v$ and $u$ are in the same cluster
of $\omega$.  There is some $a=a(G,p)<1$
such that $\tau(v,u)\le a^{ d(v,u)}$
for every $v,u\in\verts(G)$.
\endprocl

\proof
There are constants
$b=b(G,p)$ and $c=c(G,p)>0$ such that if $L$ and $L'$ are
two hyperbolic lines and the minimal distance
between $L$ and $L'$ is at least $b$,
then with probability at least $c$ there is 
an infinite path in $\d \omega$ which separates
$L$ from $L'$ in $\hh$ and does not intersect $L\cup L'$.
This follows from \ref c.half/ (or rather its generalization
to the quasi-transitive setting), since,
when $b$ is large enough,
for given $L$ and $L'$, there are hyperbolic
half spaces $H$ and $H'$, disjoint from $L\cup L'$,
such that $H\cup H'$ separates $L$ from $L'$.
Moreover, the automorphisms of $G$ act co-compactly
on $\hh$, so it is enough to consider a compact collection
of pairs $(L,L')$.

If $v,u\in\verts(G)$ and $d(v,u)$ is large, then
one can find a collection of hyperbolic lines $L_0,L_1,\dots,L_m$,
with $d(v,u)\le O(1)m$,
such that the distance between $L_j$ and $L_{j+1}$ is $b$,
and for each $j=1,2,\dots,m$,
the line $L_j$ separates $\{v\}\cup L_0\cup\cdots\cup L_{j-1}$
from $L_{j+1}\cup\cdots\cup L_m\cup \{u\}$.
The corollary immediately follows, since, by independence,
with probability at least $1-(1-c)^m$ there is
some $j=1,2,\dots,m$ and a path $\gamma\subset\d\omega$
which separates $L_{j-1}$ from $L_j$, and thereby
separates $u$ from $v$.
\Qed

\bsection{Mass Transport in the Hyperbolic Plane and Some Applications}{s.mtphyp}

The Mass Transport Principle \cite{Hag:deptree}
has an important role
in the study of percolation on nonamenable transitive graphs.
See, e.g.\ \BLPSgip. 
We now develop a continuous version of this principle,
in the setting of the hyperbolic plane, and later produce several applications.

\procl d.mtp 
A {\bf diagonally-invariant} measure $\mu$ on $\hh\times\hh$
is a measure satisfying
$$
\mu(gA\times gB) = \mu(A\times B)
$$
for all measurable $A,B\subset\hh$ and $g\in\isomh$.
\endprocl

\procl t.mtph2 \pname{Mass Transport Principle in $\hh$}
Let $\mu$ be a non-negative diagonally-invariant Borel measure on
$\hh\times\hh$.
Suppose that $\mu(A\times\hh)<\infty$
for some open $A\subset\hh$.
Then 
$$
\mu(B\times\hh)=\mu(\hh\times B)
$$
for all measurable $B\subset\hh$.
Moreover, there is a constant $c$ such that
$\mu(B\times\hh)=c\,\area(B)$ for all measurable $B\subset\hh$.

The same conclusions apply in the case where
$\mu$ is a diagonally-invariant Borel signed measure on
$\hh\times\hh$,
provided that $\left|\mu\right|(K\times\hh)<\infty$ for every compact
$K\subset\hh$.
\endprocl

Recall that the subgroup ${\rm Isom}_+(\hh)$ of orientation
preserving isometries is a simple group.  
Consequently, ${\rm Isom}_+(\hh)$ must be
contained in the kernel of the modular function,
which is a homomorphism from $\isomh$ to the multiplicative group $\R_+$.
It follows that the modular function is identically $1$;
that is, $\isomh$ is unimodular.
Let $\eta$ denote Haar measure of $\isomh$.
The particular consequence of unimodularity that we shall need is that
$\eta(A) = \eta\big\{g^{-1}\st g\in A\big\}$
for every measurable $A\subset\isomh$.
(Indeed, set $\nu(A):=\eta\left\{g^{-1}\st g\in A\right\}$.
Then it is easy to see that $\nu$ is left-invariant,
and by the uniqueness of Haar measure it must be a multiple of $\eta$. 
By choosing $A$ symmetric with respect to $g\to g^{-1}$,
it follows that $\nu=\eta$.)

\proof
We first prove the non-negative case.
Note that $\nu(B)=\mu(B\times\hh)$ is an
$\isomh$-invariant Borel measure on $\hh$.
Hence $\nu=c\,\area$, for some constant $c$.
Suppose for the moment that $\mu(\hh\times B)$ is
finite for some open $B$.
Then $\mu(\hh\times B) =c'\area(B)$,
and it remains to show that $c=c'$.

Let $B$ be some open ball in $\hh$.
Fix a point $o\in\hh$, and let $\eta$
be Haar measure on $\isomh$. 
Note that $\eta\{g\st z\in gB\}=\eta\{g\st o\in gB\}$ when $z\in\hh$.
Hence, Fubini gives, 
$$
\int \mu\big((g B)\times B\big)\, d\eta(g)=
\int_{(z,g)} 1_{z\in gB}\, d\mu(\{z\}\times B)\,d\eta(g) =
\mu(\hh\times B)\,\eta\{g\st o\in gB\}\,.
$$
Therefore,
$$
\mu(\hh\times B) =
{\int \mu\big((g B)\times B\big) d\eta(g)
\over
\eta\Big\{g\st o\in gB\Big\}
}
=
{
\int \mu\big(B\times(g^{-1} B)\big) d\eta(g)
\over
\eta\Big\{g\st o\in gB\Big\}
}
=
\mu(B\times\hh) 
\,.
$$

To complete the proof for the non-negative case, we only need to show that
$\mu(\hh\times B)<\infty$ for some open $B$.
Indeed, for every $r>0$ let $F_r$ be the set of
$(x,y)\in\hh\times\hh$ such that $d(x,y)<r$.
Set $\mu_r(K) = \mu(K\cap F_r)$.
Then $\mu_r$ is diagonally invariant and the above proof applies to it.
Therefore,
$$
\mu(\hh\times B) = \sup_{r>0} \mu_r(\hh\times B)
= \sup_{r>0} \mu_r(B\times\hh) = \mu(B\times\hh)\,.
$$
This completes the proof in the non-negative case,
and the signed version easily follows by decomposing
the measure as a difference of two non-negative measures.
\Qed

Of course, there is nothing special about the hyperbolic plane in this case;
there is a similar version of the Mass Transport Principle in any symmetric space.

Also, the Mass Transport Principle holds when the assumption
that $\mu$ is invariant under the diagonal action of
$\isomh$ is replaced by the weaker assumption that
$\mu$ is invariant under the diagonal action
of ${\rm Isom}_+(\hh)$, the group of orientation preserving isometries.
Similarly, the results stated below assuming
invariance under $\isomh$ are equally valid assuming
${\rm Isom}_+(\hh)$-invariance.

We now illustrate the Mass Transport Principle with an application
involving a relation between 
the densities of vertices, edges and faces in tilings.

Let $T$ be a random tiling of $\hh$, whose law
is invariant under $\isomh$, and with the property
that a.s.\ the edges in $T$ are piecewise smooth, and 
each face and each vertex has a finite number of edges incident with it.

\procl d.dens \procname{Vertex, face, and edge densities}
How does one measure the abundance of vertices in $T$?
Given a measurable set $A\subset\hh$,
let $N_v(A)$ be the number of vertices of $T$ in $A$. 
If $\E N_v(A)$ is finite for every measurable bounded $A\subset\hh$,
then we say that $T$ has {\bf finite vertex density}.
In this case, $\mu(A):=\E N_v(A)$ is clearly an invariant Borel measure;
hence there is a constant $\Dvert $ such that $\E N_v=\Dvert \area$.
This number $\Dvert $ will be called the {\bf vertex density}
of $T$.  The definitions of the
face and edge densities are a bit more indirect.
Fix a point $o\in\hh$.
The {\bf face density} is defined by
$\Dface =\E[ A_o^{-1}]$, where $A_o$ is the area of the tile of
$T$ that contains $o$.
Given a measurable set $A\subset\hh$, let
$Y_e(A)$ be the sum of the degrees of the vertices of $T$ in $A$. 
If $\E Y_e(A)$ is finite for bounded measurable $A$,
and $\E Y_e = 2\Dedge  \area$,
then we say that $T$ has {\bf finite edge density} $\Dedge $.
It is left to the reader to convince herself or himself that
these definitions are reasonable.
\endprocl

Given a set $B\subset\hh$ with reasonably smooth boundary 
(piecewise $C^2$ is enough),
let $\kappa_{\bd B}$ denote the associated curvature measure.
That is, for all open $A\subset \hh$,
$\kappa_{\bd B}(A)$ is the (signed) curvature of $A\cap\partial B$.
Note that if $p$ is a point  on $\bd B$ and the internal angle
of $B$ at $p$ is $\alpha$, then $\kappa_{\bd B}(\{p\})=\pi-\alpha$.
The following instance of the Gauss Bonnet Theorem will be
very useful for us.

\procl t.gbh \pname{Gauss-Bonnet in $\hh$}
Let $A$ be a closed topological disk in $\hh$,
whose boundary is a piecewise smooth simple closed path.
Then
$$
2\pi + \area(A)  =  \kappa_{\bd A}(\bd A)
\,.
$$
\endprocl

Suppose that $K\subset\hh$ is a compact set whose boundary
$\partial K$ is a $1$-manifold. If $\bd K$ has
finite unsigned curvature, namely $|\kappa_{\bd K}|(\hh)<\infty$,
then set
$$
\mu_K(A\times B) := {\area(A\cap K)\over\area(K) } \kappa_{\bd K}(B)
\,,
\label e.mukdef
$$
for every measurable $A,B\subset\hh$.
This is a signed measure on $\hh\times\hh$, and
will be very useful below.
Note that $\mu_K(\hh\times B)=\kappa_{\bd K}(B)$.

Given a tiling $T$ of $\hh$, let $\faces(T)$ denote the set of faces (that is, tiles) of $T$,
and set
$$
\mu_T = \sum_{f\in \faces(T)} \mu_f
\,, \qquad
\hat{\mu}_T = \sum_{f\in \faces(T)} |\mu_f|
\,.
$$
If $T$ is a random tiling such that
$$
\E\hat\mu_T(\hh\times A)<\infty
$$
for bounded open sets $A$, then we say that
$T$ has {\bf locally integrable curvature}.

\procl t.densities \procname{Euler formula for random tilings in $\hh$}
Let $T$ be a random tiling of $\hh$,  whose distribution
is $\isomh$-invariant.
Suppose that a.s.\ each face of $T$ is a closed topological disk
with piecewise smooth boundary, and each vertex has degree at least $3$.
Further suppose that $T$ has locally integrable curvature.
Then $T$ has finite vertex, edge and face densities,
and these densities satisfy the relation
$$
2\pi (\Dface -\Dedge +\Dvert ) = -1
\,.
\label e.densityrel
$$
\endprocl

\proof
Suppose that $\gamma$ is an arc on the boundary of two tiles,
$f$ and $f'$ of $T$, and that $\gamma$ is disjoint from the vertices
of $T$.  Then $\kappa_f(\gamma)=-\kappa_{f'}(\gamma)$, as the curvature negates under
a change of orientation.
Consequently, if $A\subset\hh$ is disjoint from the vertices of $T$,
then $\mu_T(\hh\times A)=0$.  On the other hand,
if $v$ is a vertex of a face $f$, then
$\mu_f(\hh\times\{v\})$ is equal to $\pi$ minus the
interior angle of $f$ at $v$. 
Consequently, $\mu_T(\hh\times\{v\})=\pi(\deg_v-2)$.
Since $T$ has locally integrable curvature and every vertex
has degree at least $3$,
this shows that $T$ has finite vertex and edge density, and that
$$
\E\mu_T(\hh\times A) = 2 \pi (\Dedge -\Dvert )\area(A)
\,.
$$
By the Mass Transport Principle,
$$
\E\mu_T(A\times\hh) = 2 \pi (\Dedge -\Dvert )\area(A)
\,.
$$
On the other hand, Gauss-Bonnet shows that
when $A$ is contained in a face $f\in\faces(T)$, we have
$\mu_T(A\times\hh) = \area(A)(2\pi+\area(f))/\area(f)$.
Consequently,
$$
\area(A)(2\pi \Dface +1)
=
\E\mu_T(A\times\hh) = 2 \pi (\Dedge -\Dvert )\area(A)
\,.
\Qed
$$

\procl r.sr
As the proof shows, the $-1$ on the right side
of \ref  e.densityrel/ comes from the
curvature of $\hh$.  An analogous equation holds in $\R^2$ and
in $S^2$, with the $-1$ changed to $0$ and $1$, respectively.
(For $S^2$ this is just the classical Euler formula.)
A similar proof applies to $S^2$ and $\R^2$, but for these
spaces one can also replace the use of the Mass Transport
Principle with amenability.
\endprocl

\procl r.triang
When a.s.\ all the vertices in $T$ have degree $3$, we have
$2\Dedge = 3 \Dvert$, and therefore \ref e.densityrel/ simplifies
to
$$
\Dvert =2\Dface +{1\over\pi}
\,.
\label e.dtriang
$$
\endprocl

\bsection{Voronoi Percolation in the Hyperbolic Plane}{s.hyp}

We now describe a very natural continuous percolation process 
in the hyperbolic plane.

Given a discrete nonempty set of points $X\subset\hh$,
the associated {\bf Voronoi tiling} of $\hh$
is defined by $T=T(X)=\{T_x\st x\in X\}$, where
$$
T_x:= \bigl\{y\in\hh\st \dist(y,x)=\dist(y,X)\bigr\}\,.
$$
The point $x\in X$ is called the nucleus of $T_x$.
Fix some parameters $p\in[0,1]$ and $\lambda\ge 0$.
Let $W$ be a Poisson point process in $\hh$ with intensity
$\lambda_W:=p\lambda$, and let $B$ be an independent Poisson point process
with intensity $\lambda_B:=(1-p)\lambda$.  Note that $X:=W\cup B$ is
a Poisson point process with intensity $\lambda$,
and given $X$, each point $x\in X$ is in $W$ with probability $p$,
independently.  Let $T=T(X)$ be the Voronoi tiling associated with $X$, and
set
$$
\hat W:=\bigcup_{w\in W} T_w,\qquad
\hat B:=\bigcup_{b\in B} T_b\,.
$$
Observe that a.s.\ each vertex of the tiling $T(X)$ has degree $3$.
A.s.\ $\hh$ is the union of $\hat W$ and $\hat B$, and 
$\hat W\cap\hat B$ is a $1$-manifold.  This model will be
referred to as $(p,\lambda)$-{\bf Poisson-Voronoi-Bernoulli} percolation in
$\hh$, or just Voronoi percolation, for short.
The connected components of $\hat W$ and of $\hat B$ will be called
{\bf clusters}.

It is clear that if $\hat W$ has infinite components
with positive probability, then it has infinite components
a.s.\ 
Set
\begineqalno
\Cal W
&:= \bigl\{(p,\lambda)\in[0,1]\times(0,\infty)\st
\hat W\hbox{ has infinite components a.s.}\bigr\}
\,,\cr
\Cal B
&:= \bigl\{(p,\lambda)\in[0,1]\times(0,\infty)\st
\hat B\hbox{ has infinite components a.s.}\bigr\}
\,.
\endeqalno
For every $\lambda>0$ define
$$
p_c(\lambda):= \inf\bigl\{p\in[0,1]\st (p,\lambda)\in\Cal W\bigr\}\,.
$$
It is clear that $(p,\lambda)\in\Cal W$ if $p>p_c(\lambda)$.

\proofof t.phases
It is easy to adapt the proof of \ref t.numsBern/ to this
setting.  The details are left to the reader.
\Qed

As we shall see, $[0,1]\times(0,\infty)=\Cal W\cup\Cal B$,
so \ref t.phases/ covers all possibilities.

It is easy to see that in $\Cal W\cap\Cal B$ a.s.\
all unbounded components of $\hat W$ have a cantor set
of limit points in $\hbd$, with dimension smaller than $1$.

\procl r.dfvor
One can easily show that the face density $\Dface$ of $T$ is $\lambda$,
using the following mass transport.
For each tile $f$ of $T$ with nucleus
$x$ let $\nu_f(A\times A') = \area(A\cap f)/\area(f)$ if
$x\in A'$ and $0$ otherwise.  Set $\nu_T=\sum_f\nu_f$.
Then the Mass Transport Principle with $\E\nu_T$
gives $\lambda=\Dface$.
\endprocl

\procl t.pcvor
$(1/2,\lambda)\in\Cal W\cap\Cal B$ for every $\lambda\in(0,\infty)$.
\endprocl

We present two proofs of this theorem, one uses the Mass Transport
Principle, and the other uses hyperbolic surfaces.
We start with the latter.

\medbreak\noindent{\it Hyperbolic Surfaces Proof.\enspace}
Fix some large $d>0$.  Let $S$ be a compact hyperbolic surface,
whose injectivity radius is greater than $5 d$; that is,
any disk of radius $5 d$ in $S$ is isometric to a disk in the
hyperbolic plane.  It is well known that such surfaces
exist. (See, e.g., Prop.~1 and Lemma~2 from \ref b.SSp/.)
 
Consider $(1/2,\lambda)$ percolation in $S$.
Let $K$ be the union of all white or black clusters of
diameter less than $d$.
We claim that each component of $K$ has diameter less than $d$.
Indeed, if $A$ is a white or black cluster with diameter less than $d$,
then $A$ is contained in a disk in $S$ which is isometric to a disk in
the hyperbolic plane.  It follows that the complement of
$A$ consists of one component of diameter greater than $d$,
and possibly several components of diameter smaller than $d$.
So, if a black and a white cluster have diameters $<d$ and
are adjacent, then one of them `surrounds' the other,
in the sense that the latter is contained in a component
of the complement of the first which has diameter $<d$.
It follows that each component of $K$ indeed has diameter $<d$.
 
For all $t>0$ let $K_t$ be the set of
points in $K$ with distance at least $t$ to
$S-K$.  Because each component of $K$ is isometric to a set in the
hyperbolic plane, the linear isoperimetric inequality for
the hyperbolic plane implies that
$$
\operatorname{length}(\partial K_\delta) \geq 
c \operatorname{area} K_\delta
$$
holds for all $\delta\geq 0$, where $c>0$ is some fixed
constant.
Consequently,
\begineqalno
\operatorname{area}(S)-\operatorname{area}(K_1)
&
\geq \operatorname{area}(K-K_1)
= \int_0^1\left | {d\over dt} \operatorname{area}(K_t)\right|\,dt
\cr
&= \int_0^1\operatorname{length} \partial K_t\,dt
\geq c\int_0^1 \operatorname{area}(K_t)\,dt
\geq c\operatorname{area}(K_1),
\endeqalno
which implies,
$$
\operatorname{area}(K_1) \leq (1+c)^{-1} \operatorname{area}(S).
$$
Therefore, a uniform-random point in $S$ has probability at
least $1-(1+c)^{-1}$ to be within distance $1$ of a cluster
with diameter $\geq d$.
Because the injectivity radius of $S$ is $\geq 5d$,
the same would be true for an arbitrary point in the hyperbolic plane.
Letting $d\to\infty$, we see that for any fixed point
in $\hh$, the probability that it is within distance
$1$ of an unbounded cluster is at least $1-(1+c)^{-1}$.
This implies that $(1/2,\lambda)\in\Cal W\cup \Cal B$,
and hence $(1/2,\lambda)\in\Cal W\cap \Cal B$, by symmetry.
\Qed

\medbreak\noindent{\it Mass Transport Proof.\enspace}
Consider Voronoi percolation with parameters $(1/2,\lambda)$.
Let $d>0$ be large, and let $F_d$ be the union of all
black or white components of diameter less than $d$.
Then each component $K$ of $F_d$ is a.s.\ a topological
disk with diameter bounded by $d$. 

Given a component $K$ of $F_d$, as above, let
$\mu_K$ be the signed measure on $\hh\times\hh$ defined by
$$
\mu_K(A\times B) = {\area(A\cap K)\over\area(K)} \kappa_{\partial K}(B)
\,,
$$
where $\kappa_{\partial K}$ is the curvature measure on $\partial K$.
Given the percolation configuration  $X$, let
$$
\mu^X := \sum_K \mu_K,\qquad \mu := \E\mu^X
\,,
$$
where the sum extends over all components $K$ of $F_d$.
Note that $|\mu^X(\hh\times A)|$ is bounded by
$2\pi$ times the number of vertices of the Voronoi
tiling that are in the intersection of $A$ with the boundary of $F_d$.
Consequently, $\left|\mu\right|(\hh\times A)$ is finite for every bounded $A$.
The measure $\mu$ is clearly  invariant under the diagonal
action of $\isomh$ on $\hh\times\hh$.
Therefore, the Mass Transport Principle applies to $\mu$.

Fix some point $o\in\hh$.
Fubini and the  Gauss Bonnet \ref t.gbh/ show that
$ \mu(A\times\hh)\geq \area(A) \P[o\in F_d] . $
By the Mass Transport Principle, we therefore have,
$$
\mu(\hh\times A)\geq \area(A) \P[o\in F_d] 
\,. 
\label e.muar
$$
Fix some $A$ with $\area(A)>0$.  Let $F_\infty=\bigcup_{d>0}F_d$.
{}From \ref e.muar/ we find that
$$
2\pi\E\Bigl|\bigl\{\hbox{Voronoi vertices in }A\cap \partial F_d\bigr\}\Bigr|
\geq \area(A) \P[o\in F_d] 
\,.
$$
Because $\E\bigl|\{\hbox{Voronoi vertices in }A\}\bigr|<\infty$,
by letting $d\to\infty$ it follows that
$$
2\pi\E\Bigl|\bigl\{\hbox{Voronoi vertices in }A\cap
\partial F_\infty\bigr\}\Bigr|
\geq \area(A) \P[o\in F_\infty] 
\,.
$$
This shows that $\partial F_\infty$ is not empty with positive probability.
On this event, $\hat W$ has an unbounded component or $\hat B$ has
an unbounded component.  This gives
$(1/2,\lambda)\in\Cal W\cup\Cal B$.  Symmetry then implies
that $(1/2,\lambda)\in\Cal W\cap\Cal B$, which completes the proof.
\Qed

Note that the latter proof does not require that
the tile colors be independent. 
In fact, it gives the following generalization.

\procl t.gencont
Suppose that $Z\subset\hh$ is a closed random subset
whose distribution is $\isomh$-invariant,
such that $\partial Z$ is a.s.\ a $1$-manifold
and $\Ebig{|\kappa_{\partial Z}|(A)}<\infty$ for some nonempty open
$A\subset\hh$.  Then a.s.\  there is an infinite
component in $Z$ or in $\hh\setminus Z$.
\Qed
\endprocl

We now work a bit harder to get the following explicit upper bound on $p_c(\lambda)$.
The bound will also show that $p_c(\lambda)\to 0$ as $\lambda\searrow 0$.

\procl t.vorpcbd 
$$
\forall\lambda>0,\qquad
p_c(\lambda)\leq 
{1\over 2} - {1\over 4\lambda\pi +2}
\,.
$$
\endprocl

Note that the right hand
side is zero at $\lambda=0$ and its derivative there is $\pi$.

\proof
Let $\lambda>0$, and let $p<p_c(\lambda)$.  Then $p<1/2$, by \ref t.pcvor/.
Consider three independent Poisson point processes $W,B,R$ in $\hh$,
with (positive) intensities $p\lambda,p\lambda$ and $(1-2p)\lambda$,
respectively.
Then the intensity of $W\cup B\cup R$ is $\lambda$.
Let $X=(W, B, R)$, and consider the Voronoi tiling $T=T(W\cup B\cup R)$. 
Color the tiles with white, black or red, depending on whether the
nucleus is in $W,B$ or $R$, respectively.
Let 
$$
p_R:=1-2p
\,,\qquad
p_W:=p
\,,\qquad
p_B:=p
\,.
$$
Then, given $W\cup B\cup R$,
$p_R$, $p_W$ and $p_B$ are the probabilities that any given tile
of $T(W\cup B\cup R)$ is red, white and black, respectively.

Because $p<p_c$, a.s.\ there are no unbounded white or black clusters. 
Given a cluster $K$, let the {\bf hull} of $K$ be the union
of $K$ with the bounded components of $\hh-K$.
Call a white or black cluster $K$ an {\bf empire},
if there is no black or white cluster $K'$ such that
the hull of $K'$ contains $K$. 
Let $\Cal K$ be the set of all hulls of empires.

Condition on the triplet $X=(W,B,R)$, and let $K\in\Cal K$.
Let $\kappa_{\partial K}$ denote the curvature measure on
$\partial K$.
Let $\mu_K$ be the signed measure on $\hh\times\hh$
defined by 
$$
\mu_K(A\times A') = {\area(A\cap K)\over \area(K)}
\kappa_{\partial K}(A')
\,,
$$
and let
$$
\mu^X = \sum_{K\in\Cal K} \mu_K\,,
\qquad
\mu = \E \mu^X
\,,
$$
as above.
As before, it is not hard to verify that
$|\mu|(\hh\times A)$ is finite for every bounded measurable $A$.

Fix a point $o\in\hh$. 
The Gauss Bonnet \ref t.gbh/ shows
that 
$$
\mu(A\times\hh)\geq \P[o\in\cup\Cal K] \,\area(A)\,.
\label e.aa
$$

Note that the measure $\mu^X(\hh\times \cdot)$ on $\hh$
is supported on the vertices of the Voronoi tiling
that are on the outer boundaries of the empires. 
Let $z$ be such a vertex.
Note that if $z$ does not belong to a red tile,
then $z$ is in the interior of the union of
two empires.  In that case,
$\mu^X(\hh\times\{z\})=0$, since the contributions
to $\mu^X(\hh\times\{z\})$ from both empires cancel.
Consequently, $\mu^X(\hh\times \cdot)$ is supported
on the vertices that are on the boundaries of red tiles.
Suppose that $v$ is a Voronoi vertex 
and there are three tiles $T_1$, $T_2$, $T_3$ meeting at $v$,
and $\alpha_1,\alpha_2$ and $\alpha_3$ are their angles
at $v$, respectively.  Since the Voronoi tiles are convex,
we have $\alpha_1,\alpha_2,\alpha_3\in (0,\pi]$.

If  $T_1$ is red, $T_2$ is white, and $T_3$ is black,
then $\mu^X(\hh\times\{v\}) =\pi-\alpha_2+\pi-\alpha_3=\alpha_1$
or $\mu^X(\hh\times\{v\})=0$ (the latter happens if
$v\notin\partial\cup \Cal K$).
If  $T_1$ is red, $T_2$ and $T_3$ are both of the same color, then 
$\mu^X(\hh\times\{v\})=\pi-\alpha_2-\alpha_3=\alpha_1-\pi\le 0$
or $\mu^X(\hh\times\{v\})=0$.
If $T_1$ and $T_2$ are red,
then $\mu^X(\hh\times\{v\})=\pi-\alpha_3=\alpha_1+\alpha_2-\pi$
or $\mu^X(\hh\times\{v\})=0$.
If there is no red tile among $T_1,T_2,T_3$ or if all
of them are red, then $\mu^X(\hh\times\{v\})=0$.
Consequently, given the tiling (but not the colors), the expected value of
$\mu^X(\hh\times\{v\})$ is bounded by 
$$
2p_Rp_Bp_W(\alpha_1+\alpha_2+\alpha_3)
+p_R^2(1-p_R)(2\alpha_1+2\alpha_2+2\alpha_3-3\pi)
= 2\pi p (1-2p)
\,.
$$
Hence,
$$
\mu(\hh\times A)\le \Dvert 2\pi p(1-2p)\,\area(A)
\,.
$$
Using the Mass Transport Principle and \ref e.aa/, this gives
$$
\P[o\in\cup\Cal K] \le \Dvert 2\pi p(1-2p)
\,.
\label e.lambs
$$

We now show 
$$
\P[o\in \cup\Cal K] \geq p_B+p_W=2p
\,.
\label e.empire
$$
For this, we prove that every white or black tile is
contained in the hull of an empire.  If not, there is a sequence
of black or white clusters $K_1,K_2,\dots$ such that each
$K_{j}$ is contained in the hull of $K_{j+1}$.
If infinitely many of these clusters are white, say,
then when all the red tiles are changed to black, there is still no
unbounded black cluster.  However, $p_R+p_B>1/2\ge p_c(\lambda)$,
which is a contradiction.  A similar contradiction is obtained
if infinitely many of the clusters $K_j$ are black.
Hence \ref e.empire/ holds.

Since $\Dface=\lambda$, by \ref r.dfvor/,
combining \ref e.lambs/, \ref e.empire/ and \ref e.dtriang/ gives
$$
(2\lambda\pi +1) (1-2p) \geq 1\,.
$$
This inequality must be satisfied for every $p<p_c(\lambda)$,
which proves the theorem.
\Qed


\procl l.pospc
$$\forall \lambda>0,\qquad p_c(\lambda)>0\,.$$
\endprocl

\proof
Let $Q$ be a set of points in the hyperbolic plane which
is maximal with the property that the distance between any two points
in $Q$ is at least $1$.  Then every open ball of radius $1$ contains
a point in $Q$.  

Consider Voronoi percolation with some parameters $(p,\lambda)$.
As always, let $W$ and $B$ denote the Poisson point processes
of the white and black nuclei, respectively.
Fix some large $R>0$.  Given
$y\in\hh$, consider the Voronoi tiling $T(y)$ with nuclei
$W\cup\{y\}\cup B$, and let $\hat W(y)$ denote the union of the
tiles of $T(y)$ with nuclei in $W\cup\{y\}$.  Let $\ev A(y)$
be the event that some component of $\hat W(y)$ intersects
the hyperbolic circle of radius $R$ with center $y$
and the hyperbolic circle of radius $1$ with center $y$.
The reason for introducing $T(y)$ and $\hat W(y)$ is that
the events $\ev A(y)$ and $\ev A(y')$ are independent
when $\dist(y,y')>4R$, because $\ev A(y)$ depends
only on the intersections of $B$ and $W$ with the closed ball
of radius $2R$ about $y$.  The analogous property does
not hold for $\hat W$ in place of $\hat W(y)$.
(That is, the event that there is a component of $\hat W$
which intersects the circles of radii $1$ and $R$ about
$y$ is not independent from the corresponding event for $y'$,
even if the distance between $y$ and $y'$ is large.)

Let $\alpha(R)$ be the probability that
the tile $S$ with nucleus $y$ in $T(y)$ intersects the circle of
radius $R$ about $y$.  We now estimate $\alpha(R)$ from
above.  Indeed, if $S$ intersects this circle at a point $z$,
then the open ball of radius $R$ about $z$ does not
contain any point in $B\cup W$.  Then there is a point
$x\in Q$ with $\dist(x,z)< 1$ such that the ball $B(x,R-1)$
does not intersect $B\cup W$.  For a given
$x$, the probability for that is $\exp\bigl(-\lambda\,\area(B(x,R-1))\bigr)$,
which is less than $\exp(-c\lambda e^{R})$ for some constant $c>0$
and all $R$ sufficiently large (recall that $\area(B(x,R))/e^R$
tends to a positive constant as $R\to\infty$).  Consequently,
\begineqalno
\alpha(R) 
& \le |B(y,R+1)\cap Q| \exp(-c\lambda e^R)
\le O(1) \area\bigl(B(y,R+1)\bigr) \exp(-c\lambda e^R)
\cr &
\le O(1) \exp(-c' \lambda e^R)\,,
\label e.alpha
\endeqalno
for large $R$.
Clearly, $\P[\ev A(y)]\ge\alpha(R)$.  However, there is
some very small $p(R)>0$ such that
$\P[\ev A(y)]\le 2\alpha(R)$ if $p<p(R)$. 

Assume that $p<p(R)$. Let $o\in\hh$ be some basepoint,
and let $\hat W_o$ be the component of $o$ in $\hat W$
(set $\hat W_o=\emptyset$ if $o\notin\hat W$).
Assuming that
$\hat W_o$ is unbounded, there is a.s.\ a path
$\gamma:[0,\infty)\to\hat W_o$ 
starting at $\gamma(0)=o$ and with $\dist(\gamma(t),o)\to\infty$ as
$t\to\infty$. 
In this case, let $t_0:=0$, $y_0:=o$, and inductively set
$t_n:=\sup\{t\st \dist(\gamma(t),y_{n-1})=5R\}$,
and $y_n:=\gamma(t_n)$.  Let $y_n'$ be a point in $Q$ closest to $y_n$.
Then for each $n$, $\dist(y_n',y_{n+1}')\le 5R+2$,
and for each $j\neq k$, $\dist(y_k',y_j')\ge 5R-2$.
Observe that the events $\ev A(y'_k)$ all hold.
Consequently, for every $n=1,2,\dots$, we may bound the
event that $o$ is in an unbounded component of $\hat W$ by
$$
\sum_{x_0,x_1,\dots,x_n} \P\bigl[\ev A(x_0)\cap\ev A(x_1)\cap\cdots \cap \ev A(x_n)\bigr]\,,
$$
where the sum is over all sequences $x_0,x_1,\dots,x_n$ in $Q$
such that $x_0$ is within distance $1$ of $o$, each $x_j$
is within distance $5R+2$ of $x_{j-1}$, and
$\dist(x_j,x_k)\ge 5R-2$ when $j\ne k$.  Assuming $R>2$,
these events 
$\ev A(x_0),\ev A(x_1),\dots , \ev A(x_n)$ are independent,
and we get the bound
$$
\Pbig{\hat W_o\hbox{ is unbounded}}
\le 
\bigl(2\alpha(R)\bigr)^{n+1} \bigl|\{\hbox{such sequences }x_0,\dots,x_n\}\bigr|\,.
\label e.aaa
$$
Now, the number of such sequences is at most 
$$
\max\Bigl\{ \bigl|Q\cap B(z,5R+2)\bigr|^n \st z\in\hh\Bigr\}
$$
times the number of possible choices of $x_0$.
This is at most $\exp(c_1 R n)$, for some constant $c_1$.
By our estimate \ref e.alpha/ for $\alpha(R)$, it is clear that there is
some large $R_0>0$ such that the right hand side
of \ref e.aaa/ goes to zero as $n\to\infty$.
Then for $p<p(R_0)$ the probability that
$\hat W_o$ is unbounded is zero, and hence
$p_c(\lambda)\ge p(R_0)$.
\Qed

\procl l.contpc
$p_c(\lambda)$ is continuous on $(0,\infty)$.
\endprocl

\proof
Note that given $\lambda,h>0$, a union of two independent Poisson point
processes with intensities $\lambda$ and $h$ is a Poisson point process
with intensity $\lambda+h$.  Fix some $o\in\hh$, and let
$\theta$ be the probability that $o$ is in an unbounded component of $\hat W$.
Consider $\theta=\theta(\lambda_W,\lambda_B)$
as a function of $\lambda_W=p\lambda$ and $\lambda_B=(1-p)\lambda$.
It is clear that $\theta$ is monotone increasing (weakly) in $\lambda_W$ and monotone
decreasing (weakly) in $\lambda_B$.  Consequently, if $\lambda'>\lambda$,
we must have $p_c(\lambda')\lambda'\ge p_c(\lambda)\lambda$ and
$\bigl(1-p_c(\lambda')\bigr)\lambda'\ge \bigl(1-p_c(\lambda)\bigr)\lambda$.
Hence,
$$
p_c(\lambda){\lambda\over\lambda'}\le p_c(\lambda')
\le 1-\bigl(1-p_c(\lambda)\bigr){\lambda\over\lambda'}\,,
$$
which implies continuity.
\Qed

\proofof t.hber
Part (a) follows from \ref t.vorpcbd/ and
\ref l.pospc/.  Part (b) follows from (a).
Part (c) is \ref l.contpc/, above.
The proof from \BLPScrit{} of
\ref t.critperc/
can easily be adapted to prove (d).
\Qed

%
%
%
%
%

\bsection{Open problems}{s.open}

\procl q.noplan
In the absence of planarity, most of the proofs in this
paper are invalid.  Which results can be extended to
transitive graphs that are quasi-isometric to the hyperbolic
plane, but are not planar?
\endprocl

\procl g.vorpc  $\lim_{\lambda\to\infty}p_c(\lambda)=1/2$.
\endprocl

\procl q.vorpc
What are the asymptotics of $p_c(\lambda)$ as
$\lambda\to\infty$ and as $\lambda\to 0$?
In particular, what is $p_c'(0)$?
\endprocl

\procl g.mon
$p_c(\lambda)$ is strictly monotone increasing.
\endprocl

Given $a>0$, let $a\hh$ denote $\hh$ with the metric scaled
by $a$.  It is clear that $(p,\lambda)$-Voronoi percolation
on $\hh$ is the same as $(p,1)$-Voronoi percolation on
${\sqrt{\lambda}}\H^2$.  As $a\to\infty$, the space
$a\hh$ looks more and more like $\R^2$ ($a\hh$
has constant curvature $-1/a^2$).  This leads to the following
question.

\procl q.euc
It is known that for Voronoi percolation in
the Euclidean plane $p_c\ge 1/2$ \cite{Zvavitch:vor} (in the
Euclidean setting, $p_c$ clearly does not depend on $\lambda$).
However, the conjecture $p_c=1/2$ is still open.
Lacking is a proof that there are unbounded components at $p>1/2$.  
Is it possible to prove that $p_c=1/2$ in the Euclidean setting,
by taking a limit as $\lambda\to\infty$ in the hyperbolic setting?
\endprocl

This direction can also lead to a plausible guess as to the
asymptotics of $p_c(\lambda)$ as $\lambda\to\infty$.
Recall the notion of the correlation length $\xi(p)$
(see \cite{Grimmett:percolation}), which roughly measures
the length scale at which Bernoulli($p$) percolation 
is substantially different from Bernoulli($p_c$) percolation.
Similarly, at lengths on the order of $a$ or higher, the space
$a\hh$ appears substantially different from $\R^2$.
On balls of size smaller than $a$, $(p,1)$-Voronoi percolation
on $a\hh$ does not look too different from
$(p,1)$-Voronoi percolation on $\R^2$.
This suggests that if $\xi(p)$ is significantly smaller
than $a$ (here $\xi(p)$ is the
correlation length for $(p,1)$-Voronoi percolation on $\R^2$),
then there will be no unbounded white clusters for $(p,1)$-Voronoi percolation
on $a\hh$.  Conversely, if $\xi(p)$ is significantly larger
than $a$, the hyperbolicity of $a\hh$ appearing
on the scale of $a$ will help to create unbounded white clusters.
Since $(p,1)$-Voronoi percolation on ${\sqrt{\lambda}}\H^2$
is the same as $(p,\lambda)$-Voronoi percolation on $\hh$,
this suggests that $\xi\bigl(p_c(\lambda)\bigr)$ should grow roughly
at the same rate as $\sqrt{\lambda}$ when $\lambda\to\infty$.
It is conjectured that $\xi(p)$ grows like $|p-p_c|^{-4/3}$;
which leads to the guess that $p_c(\lambda)$
is asymptotic to ${1\over 2} - \lambda^{-2/3}$.


\bibsty{mybibstyle}
\bibfile{\jobname}
\startbib{BLPS99b}
\begingroup

\bibitem[Bab97]{Babai:growth}
{\sc L.~Babai}.
\newblock The growth rate of vertex-transitive planar graphs.
\newblock In {\em Proceedings of the Eighth Annual ACM-SIAM Symposium on
  Discrete Algorithms (New Orleans, LA, 1997)}, pages 564--573, New York, 1997.
  ACM.

\bibitem[BK89]{BK:uni}
{\sc R.~M. Burton and M.~Keane}.
\newblock Density and uniqueness in percolation.
\newblock {\em Comm. Math. Phys. 121}, 3 (1989), pages 501--505.

\bibitem[BK91]{Burton-Keane:planar}
{\sc R.~M. Burton and M.~Keane}.
\newblock Topological and metric properties of infinite clusters in stationary
  two-dimensional site percolation.
\newblock {\em Israel J. Math. 76}, 3 (1991), pages 299--316.

\bibitem[BLPS99a]{BLPSgip}
{\sc I.~Benjamini, R.~Lyons, Y.~Peres, and O.~Schramm}.
\newblock Group-invariant percolation on graphs.
\newblock {\em Geom. Funct. Anal. 9}, 1 (1999), pages 29--66.

\bibitem[BLPS99b]{BLPSdeath}
{\sc I.~Benjamini, R.~Lyons, Y.~Peres, and O.~Schramm}.
\newblock Critical percolation on any nonamenable group has no infinite
  clusters.
\newblock {\em Ann. Probab. 27}, 3 (1999), pages 1347--1356.

\bibitem[BLPS00]{BLPS:usf}
{\sc I.~Benjamini, R.~Lyons, Y.~Peres, and O.~Schramm}.
\newblock Uniform spanning forests.
\newblock {\em Ann. Probab.\/} (2000).
\newblock To appear.

\bibitem[BLS99]{BLS:pert}
{\sc I.~Benjamini, R.~Lyons, and O.~Schramm}.
\newblock Percolation perturbations in potential theory and random walks.
\newblock In M.~Picardello and W.~Woess, editors, {\em Random Walks and
  Discrete Potential Theory}, Sympos. Math., pages 56--84, Cambridge, 1999.
  Cambridge Univ. Press.
\newblock Papers from the workshop held in Cortona, 1997.

\bibitem[BS96]{pyond}
{\sc I.~Benjamini and O.~Schramm}.
\newblock Percolation beyond ${\Z}\sp d$, many questions and a few answers.
\newblock {\em Electron. Comm. Probab. 1\/} (1996), pages no.\ 8, 71--82
  (electronic).

\bibitem[BS99]{pyond:recent}
{\sc I.~Benjamini and O.~Schramm}.
\newblock Recent progress on percolation beyond ${\Z}\sp d$. 
\hfill\break
\newblock \htmlref{http://www.wisdom.weizmann.ac.il/\string~schramm/papers/pyond-rep/}
  {{\tt http://www.wisdom.weizmann.ac.il/\string~schramm/papers/pyond-rep/}},
  1999.

\bibitem[BSt90]{BeSt}
{\sc A.~F. Beardon and K.~Stephenson}.
\newblock The uniformization theorem for circle packings.
\newblock {\em Indiana Univ. Math. J. 39}, 4 (1990), pages 1383--1425.


\bibitem[CFKP97]{CFKW}
{\sc J.~W. Cannon, W.~J. Floyd, R.~Kenyon, and W.~R. Parry}.
\newblock Hyperbolic geometry.
\newblock In {\em Flavors of geometry}, pages 59--115. Cambridge Univ. Press,
  Cambridge, 1997.

\bibitem[DS84]{DoyleSnell}
{\sc P.~G. Doyle and J.~L. Snell}.
\newblock {\em Random Walks and Electric Networks}.
\newblock Mathematical Association of America, Washington, DC, 1984.

\bibitem[GN90]{GN:treeZ}
{\sc G.~R. Grimmett and C.~M. Newman}.
\newblock Percolation in $\infty+1$ dimensions.
\newblock In G.~R. Grimmett and D.~J.~A. Welsh, editors, {\em Disorder in
  Physical Systems}, pages 167--190. Oxford Univ. Press, New York, 1990.

\bibitem[Gri89]{Grimmett:percolation}
{\sc G.~Grimmett}.
\newblock {\em Percolation}.
\newblock Springer-Verlag, New York, 1989.

\bibitem[H{\"a}g97]{Hag:deptree}
{\sc O.~H{\"a}ggstr{\"o}m}.
\newblock Infinite clusters in dependent automorphism invariant percolation on
  trees.
\newblock {\em Ann. Probab. 25}, 3 (1997), pages 1423--1436.

\bibitem[HP99]{HP:unimon}
{\sc O.~H{\"a}ggstr{\"o}m and Y.~Peres}.
\newblock Monotonicity of uniqueness for percolation on {C}ayley graphs: all
  infinite clusters are born simultaneously.
\newblock {\em Probab. Theory Related Fields 113}, 2 (1999), pages 273--285.

\bibitem[HS95]{HS:hyp}
{\sc Z.-X. He and O.~Schramm}.
\newblock Hyperbolic and parabolic packings.
\newblock {\em Discrete Comput. Geom. 14}, 2 (1995), pages 123--149.

\bibitem[Imr75]{Imrich}
{\sc W.~Imrich}.
\newblock On {W}hitney's theorem on the unique embeddability of $3$-connected
  planar graphs.
\newblock In {\em Recent advances in graph theory (Proc. Second Czechoslovak
  Sympos., Prague, 1974)}, pages 303--306. Academia, Prague, 1975.

\bibitem[Lal98]{Lalley:Fuch}
{\sc S.~P. Lalley}.
\newblock Percolation on {F}uchsian groups.
\newblock {\em Ann. Inst. H. Poincar\'e Probab. Statist. 34}, 2 (1998), pages
  151--177.


\bibitem[Lyo00]{Lyons:percsurvey}
{\sc R.~Lyons}.
\newblock Phase transitions on nonamenable graphs.
\newblock {\em J. Math. Phys. 41}, 3 (2000), pages 1099--1126.
\newblock Probabilistic techniques in equilibrium and nonequilibrium
  statistical physics.

\bibitem[Lyo01]{Lyons:book}
{\sc R.~Lyons}.
\newblock {\em Probability on Trees and Networks}.
\newblock Cambridge University Press, 2001.
\newblock Written with the assistance of Y. Peres. In preparation. Current
  version available at \htmlref{http://php.indiana.edu/\string~rdlyons/}{{\tt
  http://php.indiana.edu/\string~rdlyons/}}.

\bibitem[Mad70]{Mader}
{\sc W.~Mader}.
\newblock \"{U}ber den {Z}usammenhang symmetrischer {G}raphen.
\newblock {\em Arch. Math. (Basel) 21\/} (1970), pages 331--336.

\bibitem[MR96]{MR:contperc}
{\sc R.~Meester and R.~Roy}.
\newblock {\em Continuum percolation}.
\newblock Cambridge University Press, Cambridge, 1996.

\bibitem[PSN00]{PakSN:uniq}
{\sc I.~Pak and T.~Smirnova-Nagnibeda}.
\newblock On non-uniqueness of percolation on nonamenable {C}ayley graphs.
\newblock {\em C. R. Acad. Sci. Paris S\'er. I Math. 330}, 6 (2000), pages
  495--500.

\bibitem[Per00]{P:nonamen}
{\sc Y.~Peres}.
\newblock Percolation on nonamenable products at the uniqueness threshold.
\newblock {\em Ann. Inst. H. Poincar\'e Probab. Statist. 36}, 3 (2000), pages
  395--406.

\bibitem[Sch99a]{S:nuni}
{\sc R.~H. Schonmann}.
\newblock Percolation in $\infty+1$ dimensions at the uniqueness threshold.
\newblock In M.~Bramson and R.~Durrett, editors, {\em Perplexing Probability
  Problems: Papers in Honor of Harry Kesten}, pages 53--67, Boston, 1999.
  Birkh{\"a}user.

\bibitem[Sch99b]{Scho:unimon}
{\sc R.~H. Schonmann}.
\newblock Stability of infinite clusters in supercritical percolation.
\newblock {\em Probab. Theory Related Fields 113}, 2 (1999), pages 287--300.

\bibitem[SS97]{SSp}
{\sc P.~Schmutz~Schaller}.
\newblock Extremal {R}iemann surfaces with a large number of systoles.
\newblock In {\em Extremal Riemann surfaces (San Francisco, CA, 1995)}, pages
  9--19. Amer. Math. Soc., Providence, RI, 1997.

\bibitem[Wat70]{Watkins}
{\sc M.~E. Watkins}.
\newblock Connectivity of transitive graphs.
\newblock {\em J. Combinatorial Theory 8\/} (1970), pages 23--29.

\bibitem[Zva96]{Zvavitch:vor}
{\sc A.~Zvavitch}.
\newblock The critical probability for {V}oronoi percolation, 1996.
\newblock MSc. thesis, Weizmann Institute of Science.

\endbib
\endreferences

\filbreak
\begingroup
\eightpoint\sc
\parindent=0pt 

\medskip
{\rm Itai Benjamini
\par Weizmann Institute of Science}
\par {\tt itai@wisdom.weizmann.ac.il}
\par {\tt http://www.wisdom.weizmann.ac.il/\string~itai/}

\medskip
{\rm Oded Schramm
\par Weizmann Institute of Science and Microsoft Research}
\par {\tt schramm@wisdom.weizmann.ac.il}
\par {\tt schramm@microsoft.com}
\par {\tt http://www.wisdom.weizmann.ac.il/\string~schramm/}

\endgroup

\bye